\documentclass[a4paper, 12pt]{article}
\usepackage{amsmath,graphicx,amssymb,amsthm}
\allowdisplaybreaks
\def\cc{{\mathbb{C}}}

\def\zz{{\bf Z}}

\def\AA{\mathfrak A}

\def\A{\mathcal A}
\def\B{\mathcal B}
\def\G{\mathcal G}
\def\K{\mathcal K}

\def\D{\mathcal D}
\def\B{\mathcal B}

\def\H{\mathcal H}
\def\K{\mathcal K}
\def\L{\mathcal L}
\def\M{\mathcal M}
\def\N{\mathcal N}
\def\R{\mathcal R}
\def\S{\mathcal S}

\def\amslatex{$\mathcal{A}\kern-.1667em\lower.5ex\hbox{$\mathcal{M}$}\kern-.125em\mathcal{S}$-\LaTeX}

\def \sp {\hspace{0.3cm}}
\def\minus{\mathop{\ominus}}
\def\tensor{\mathop{\bar\otimes}}
\newtheorem{set}{set}[section]
\newtheorem{Corollary}[set]{Corollary}

\newtheorem{Lemma}[set]{Lemma}

\newtheorem{Proposition}[set]{Proposition}

\newtheorem{Theorem}[set]{Theorem}
\newcommand{\define}{\mathrel{\hbox{$\equiv$\hskip -.90em \lower .47ex \hbox{$\leftharpoondown$}}}}
\newcommand{\enifed}{\mathrel{\hbox{$\equiv$\hskip -.90em \lower .47ex \hbox{$\rightharpoondown$}}}}

\begin{document}
\pagestyle{myheadings}
\date{}
\title {On Transitive Algebras Containing a Standard Finite von Neumann Subalgebra}
\author{Junsheng Fang \and Don Hadwin\and  Mohan Ravichandran}
 \maketitle
\begin{abstract}

Let $\M$ be a finite von Neumann algebra acting on a Hilbert space
$\H$  and $\AA$ be a transitive algebra containing $\M'$. In this
paper we prove that if $\AA$ is $2-$fold transitive, then $\AA$ is
strongly dense in $\B(\H)$.  This implies that if a transitive
algebra containing a standard finite von Neumann algebra (in the
sense of~\cite{Ha1}) is $2-$fold transitive, then $\AA$ is strongly
dense in $\B(\H)$. Non-selfadjoint algebras related to free products
of finite von Neumann algebras, e.g., $\L{\mathbb{F}_n}$ and
$(M_2(\cc), \frac{1}{2}Tr)*(M_2(\cc), \frac{1}{2}Tr)$, are studied.
Brown measures of certain operators in $(M_2(\cc),
\frac{1}{2}Tr)*(M_2(\cc), \frac{1}{2}Tr)$ are explicitly computed.
\end{abstract}

\vskip1.0cm

\centerline{{\Large\bf Introduction}} \mbox{}\\

 Let $\H$ be a complex Hilbert
space and $\B(\H)$ be the set of bounded linear operators on $\H$. A
subalgebra $\AA$ of $\B(\H)$ is \emph{transitive} if it contains the
identity operator and has no invariant (closed) subspace other than
the two trivial ones. The  \emph{transitive algebra problem}
 asks: if $\AA$ is a transitive algebra on
$\H$, is $\AA$ strongly
 dense in $\B(\H)$? This problem was implicitly contained in a
question of R. Kadison~\cite{Ka} concerning algebras whose invariant
subspaces have invariant complementary subspaces, and it was W.
Arveson~\cite{Ar} who first explicitly stated the problem, coined
the term "transitive algebra", and began an in-depth study of the
problem. Note that an affirmative answer to the transitive algebra
problem would give rise to an affirmative answer to (hyper)invariant
subspace problem. The (hyper)invariant subspace problem asks if  an
algebra generated by (the commutant of) a single bounded operator on
$\H$ can be
transitive.\\

 On the other hand, one's intuition expects that there
exists a transitive algebra which is not strongly dense in $\B(\H)$.
It is of interest, then, to know how one might strengthen the
hypothesis so as to get a provable result. The first partial
solutions of the transitive algebra problem were given by
Arveson~\cite{Ar}. Arveson proved that if $\AA$ is a transitive
algebra containing a MASA (maximal abelian von Neumann subalgebra)
of $\B(\H)$, then the strong closure of $\AA$ is $\B(\H)$. For
various generalizations of Arveson's results, we refer to~\cite{D-P,
No, NRR}.  Inspired by the invariant subspace problem affiliated
with a von Neumann algebra, we consider the following construction:
Let $\M$ be a von Neumann algebra acting on a Hilbert space and
$\M'$ be the commutant of $\M$.  Suppose $\{T_{\alpha}\}\subseteq
\M$ has no nontrivial common invariant subspace relative to $\M$,
i.e, if $E\in\M$ satisfies $ET_{\alpha}E=T_{\alpha}E$ for all
$T_{\alpha}$ then $E=0$ or $I$. It is easy to see that the algebra,
$\AA$, generated by $\{T_{\alpha}\}$ and $\M'$ is a transitive
algebra. Now we may ask that if   $\AA$ is strongly dense in
$\B(\H)$ when we choose $\M$ and $\{T_{\alpha}\}$ suitably. To make
the question
non-trivial, we first choose the ``size" of $\M$ suitably large.\\

 Recall that a von Neumann
algebra $\M$ acting on a Hilbert space $\H$ is said to be
\emph{standard} if there exists a conjugate unitary operator
$J:\H\rightarrow\H$, such that the mapping $X\rightarrow JX^*J$ is a
* anti-isomorphism from $\M$ onto $\M^{'}$.  Haagerup~\cite{Ha1}
proved that every von Neumann algebra is
*-isomorphic to a standard von Neumann algebra  on a Hilbert space.
 If a von Neumann algebra $\M$ is standard  on a Hilbert space
$\H$, then $\M'$ is  also standard on $\H$.  We may ask the
following question: If a von Neumann algebra $\M$ is standard on a
Hilbert space $\H$ and $\AA$ is a transitive algebra on $\H$ which
contains $\M'$, is  $\AA$ strongly dense in  $\B(\H)$?
 By~\cite{R-R} (Theorem 8.26), if $\M$ is a type $I_{\infty}$ factor, this
question is equivalent to Kadison's transitive algebra question. So
it is more interesting to restrict one's attention to the case where
$\M$ is a finite von Neumann algebra. Notably, if  an abelian von
Neumann algebra $\M$ is standard on a Hilbert space $\H$, then $\M$
is a MASA  of $\B(\H)$. By~\cite{Ar}, any transitive algebra which
contains $\M'=\M$ is strongly dense in $\B(\H)$ and thus the answer
to above question is affirmative. Suppose $\M$ is a finite von
Neumann algebra acting on a Hilbert space $\H$ and $\AA$ is a
transitive algebra which contains $\M'$. In section 3 of this paper,
we prove the following result: if $\AA$ is 2-fold transitive (i.e.,
for any linearly independent vectors $\xi,\eta$ in $\H$, the closure
of $\{(T\xi,T\eta):\,\, T\in\AA\}$ is $\H\oplus\H$),  then  $\AA$ is
strongly dense in $\B(\H)$. As a corollary, this implies that if
$\M$ is a standard finite von Neumann algebra  on a Hilbert space
$\H$ and $\AA$ is a $2-$fold transitive algebra containing $\M$,
then  $\AA$ is strongly dense in  $\B(\H)$. This partly answers a
question of Arveson~\cite{Ar2} (also see 10.5 of~\cite{R-R}), which
asks for a transitive algebra $\AA$ whether $2-$fold transitivity
implies that the strong closure of $\AA$ is $\B(\H)$. The proof of
our result relies on a new characterization
of $n-$fold transitivity and operator theory techniques.\\

As applications, we study some non-selfadjoint algebras related to
(reduced) free products of finite  von Neumann algebras. In section
4, we prove the following result: Let$(\M,\tau)$ be a finite von
Neumann algebra with a faithful normal trace $\tau$. Suppose
$\N\subseteq\M$ is a von Neumann subalgebra and $Z\in \M$ satisfies
the following conditions:
\begin{enumerate}
\item $Z\neq 0$ and $\tau(Z)=0$; \item $Z, Z((\N^{'}\cap\M)\minus
\cc I), ((\N^{'}\cap\M)\minus \cc I)Z$ are mutually orthogonal in
$L^2(\M,\tau)$.
\end{enumerate}
Let $\AA\subseteq\B(L^2(\M,\tau))$ be the  algebra generated by $Z,
\N$ and $\M^{'}$ (relative to $\B(L^2(\M,\tau))$). Then $\AA$ is a
transitive algebra and  strongly dense in $\B(L^2(\M,\tau))$. As a
corollary, we prove the following: Let $(\M_1,\tau_1),
(\M_2,\tau_2)$ be finite von Neumann algebras and
$\M=(\M_1,\tau_1)*(\M_2,\tau_2)$ be the reduced free product von
Neumann algebra and $\tau$ be the induced faithful normal trace on
$\M$. Suppose $\N\subseteq\M_1$ is a diffuse von Neumann subalgebra
and $Z\in \M_2$ is not a scalar. Let $\AA\subseteq\B(L^2(\M,\tau))$
be the algebra generated by $Z, \N$ and $\M^{'}$. Then $\AA$ is a
transitive algebra and strongly dense in $\B(L^2(\M,\tau))$. For
example, let $\M=\L{\mathbb{F}_n}$ ($2\leq n\leq \infty$) be the
type $II_1$ factor associated with the left regular representation
$\lambda$ of the free group $\mathbb{F}_n$ on $n$ generators and $a,
b$ be two generators of $\mathbb{F}_n$. Let $k$ be a positive
integer and $Z\in\{\lambda(a)\}''$ be a nonscalar operator.  We show
that the algebra $\AA$ on $L^2(\M,\tau)$ generated by $Z,
\lambda(b)^k,$ and $\M'$ is a transitive algebra and strongly dense
in $\B(L^2(\M,\tau))$. The proof of above results relies on
techniques of type $II_1$ factors. Freeness also plays a key role.\\

In section 5, we consider some interesting (but by no means trivial)
non-selfadjoint algebras relate to $\M=(M_2(\cc),
\frac{1}{2}Tr)*(M_2(\cc), \frac{1}{2}Tr)$(reduced free product von
Neumann algebra with respect to the normalized trace on $M_2(\cc)$).
Let $\left(E_{ij}\right)_{i,j=1, 2}$ and
$\left(F_{ij}\right)_{i,j=1, 2}$ be the matrix units in $(M_2(\cc),
\frac{1}{2}Tr)*1$ and $1*(M_2(\cc), \frac{1}{2}Tr)$, respectively.
Consider the following subalgebras of $\B(L^2(\M,\tau))$:
\begin{enumerate}
\item $\AA_1$, the algebra generated by $E_{11}, E_{12}, F_{11}$ and
$\M'$;
\item $\AA_2$, the algebra generated by $E_{11}, E_{12}, F_{12}$ and
$\M'$;
\item $\AA_3$, the algebra generated by $E_{11}, F_{12}$ and $\M'$;
\item $\AA_4$, the algebra generated by $E_{12}, F_{12}$ and $\M'$.
\end{enumerate}

In this paper, we prove that $\AA_1$ is a transitive algebra and
strongly dense in $\B(L^2(\M,\tau))$. We then prove that $\AA_3$ and
$\AA_4$ are not transitive algebras. Indeed, for any $0\leq r\leq
1$, there are invariant subspaces $E_r$ and $F_r$ of $\AA_3$ and
$\AA_4$, respectively, such that $\tau(E_r)=\tau(F_r)=r$. The main
idea to prove that $\AA_3$ and $\AA_4$ are not transitive is as
following. Note that $\AA_3\subseteq \{(E_{11}-E_{22}+F_{12})^2\}'$
and $\AA_4\subseteq \{(E_{12}+F_{12})^2\}'$. Based on the results of
Brown measures of $R-$diagonal operators computed by Haagerup and
Larson~\cite{H-L}, we explicitly compute the Brown measures of
$E_{11}-E_{22}+F_{11}$ and $E_{12}+F_{12}$ by using  techniques of
free probability theory and operator theory. Then we apply Theorem
7.1 of~\cite{H-S} on the existence of hyperinvariant subspaces of
operators which have non-single point Brown measure relative to
factors of type $II_1$.
The question that $\AA_2$ is transitive or not remains open! \\

Remarkable progress on the (hyper)invariant subspace problem
relative to a factor of type $II_1$ have been made during past ten
years (see for example~\cite{D-H1,D-H2,Ha2,H-S,S-S}).  The fact that
2-fold transitivity of a transitive algebra  containing a standard
type $II_1$ factor implies that the strong closure of the transitive
algebra is $\B(\H)$ can be viewed as a support for the existence of
nontrivial (hyper)invariant subspaces of operators relative to a
factor of type
$II_1$.\\

Besides the introduction, section 3, section 4 and section 5, there
are two more sections in this paper. In section 1, we provide some
characterizations of $n-$fold transitivity. To prove our main result
(Theorem 3.1), some auxiliary lemmas are proved in section 2. \\

For the general theory of operator theory and invariant subspaces,
we refer to~\cite{R-R}. For the general theory of von Neumann
algebras, we refer to~\cite{K-R}. For the general theory of free
probability
theory, we refer to~\cite{VDN}.\\

We wish to thank professor Liming Ge for valuable suggestions that
originally aroused our interest in the transitive algebra problem.
We also thank professor William Arveson for valuable comments on the
transitive algebra problem.

\section{On $n-$fold transitivity}
We begin this section by establishing some notation and terminology.
Let $\H$ be a  Hilbert space. A \emph{linear manifold} in $\H$ is a
subset of $\H$ which is closed under vector addition and under
multiplication by complex numbers. A \emph{subspace} of $\H$ is a
linear manifold which is closed in the norm topology; the trivial
subspaces are $\{0\}$ and $\H$. If $\D$ is a linear manifold in
$\H$, then $[\D]$ denotes the norm closure of $\D$. \\

 If $T\in\B(\H)$, the collection of
all subspaces of $\H$ invariant under $T$ is denoted by $Lat T$; if
$\AA\subseteq\B(\H)$, then $Lat\AA=\cap_{T\in\AA} Lat T$.  A
subspace $\K$ is \emph{hyperinvariant} for $T$ if $\K\in Lat\{T\}'$,
i.e., $\K\in Lat S$ for every $S$ commute with $T$. Let $P\in\B(\H)$
be a projection, i.e., $P=P^*=P^2$. We say $P\in Lat T$ if
$P\H\in Lat T$. $P\in Lat T$ if and only if $PTP=TP$. \\

If $\H$ is a Hilbert space and $n$ is a positive integer, then
$\H^{(n)}$ denotes the direct sum of $n$ copies of $\H$, i.e.,
 the Hilbert space $\H\oplus\cdots\oplus\H$. If $T$ is an operator
 on $\H$, then $T^{(n)}$ denotes the direct sum of $n$ copies of
 $T$ (regarded as an operator on $\H^{(n)}$). However, we will use $I_n$
 instead of $I^{(n)}$ to denote the identity operator on $\H^{(n)}$. If $\AA$ is a set of
 operators on $\H$, then $\AA^{(n)}=\{T^{(n)}:\,\, T\in\AA\}$. We
 will identify $\B(\H^{(n)})$ with $M_n(\mathbb{C})\bar{\otimes} \B(\H)$ by
 writing $T\in \B(\H^{(n)})$ as matrix form $(T_{ij})_{n \times n}$.
 With this identification, $\AA^{(n)}=\mathbb{C}I_n\bar{\otimes}\AA$.\\

Let $\AA\subseteq\B(\H)$ be a transitive algebra and $n\in
\mathbb{N}$. Recall that $\AA$ is said to be \emph{$n-$fold
transitive} if for any linearly independent vectors $\xi_1, \xi_2,
\cdots, \xi_n$ in $\H$, $[(S\xi_1,S\xi_2,\cdots, S\xi_n):\,\,
S\in\AA ]=\H^{(n)}.$ Note that if $\AA$ is $n-$fold transitive, then
it is also $m-$fold transitive for each $m<n$. We consider $n=2$,
first. The following lemma is well known (cf.~\cite{Ar}). For the
sake of completeness, we provide the proof.

\begin{Lemma} Let $\AA$ be a transitive algebra on a Hilbert space
$\H$. For any $\xi,\eta\in\H$, $\xi,\eta\neq 0$, either
$[(T\xi,T\eta):\,\, T\in\AA]=\H^{(2)}$ or there is a closed, densely
defined operator $S$ such that $ST=TS$ for every $T\in\AA$ and
 $\G(S)$, the graph of $S$, equals to
$[(T\xi,T\eta):\,\, T\in\AA]$. If $\G(S)=[(T\xi,T\eta):\,\,
T\in\AA]$, then $S^{-1}$ \emph{(}as a mapping\emph{)} exists and
$\G(S^{-1})=[(T\eta,T\xi):\,\, T\in\AA]$.
\end{Lemma}
\begin{proof}We can assume that $\xi, \eta$ are linearly
independent. Suppose $\G=[(T\xi,T\eta):\,\, T\in\AA]\neq\H^{(2)}$.
If there is $z\neq 0$ such that $(0,z)\in\G$, then the closure of
$\{(0,Tz):\,\, T\in\AA\}$ is $0\oplus \H$  since $\AA$ is
transitive. Thus $0\oplus \H\subseteq \G$ and $\G=\H^{(2)}$. It is a
contradiction. So $(0,z)\in\G$ implies that $z=0$.  Define
$ST\xi=T\eta$. Then $S$ is well defined and $\G(S)=\G$. So $S$ is a
closed, densely defined operator. By symmetry of $\xi$ and $\eta$,
$S^{-1}T\eta=T\xi$ is a closed, densely defined operator such that
$\G(S^{-1})=[(T\eta,T\xi):\,\, T\in\AA]$.
\end{proof}

 In the following proposition we summarize some
characterizations of $2-$fold transitivity. ``$1\Leftrightarrow 2$''
is proved by Arveson in~\cite{Ar}. The authors can not find the
equivalence of $1, 3, 4$ in the literature(for example~\cite{Ar}
and~\cite{R-R}). For the sake of completeness, we provide the proof.
\begin{Proposition}Let $\AA\subseteq\B(\H)$ be a transitive algebra. Then
the following conditions are equivalent:
\begin{enumerate}
\item $\AA$ is $2-$fold transitive;
\item If $S$ is a closed, densely defined operator such that for any
$T\in\AA$,  $TS=ST$, i.e., for any $\xi\in \D(S)$, $TS\xi=ST\xi$,
then $S=\lambda I$ for some scalar $\lambda$;
\item For any $\xi,\eta$ and $\zeta\in\H$, $\xi\neq 0$, there exist a sequence $T_n$ in $\AA$
such that $T_n\xi$ converges to $\zeta$ and $\sup_n\|T_n\eta\|<
\infty$;
\item Lat $\AA^{(2)}$=Lat $\B(\H)^{(2)}$, i.e., for any projection $P\in
M_2(\cc)\tensor \B(\H)$, if $P\in\text{Lat}(\cc I_2\tensor \AA)$,
then $P\in M_2(\cc)\tensor \cc I$.
\end{enumerate}
\end{Proposition}
\begin{proof}
``$1\Rightarrow 2$''.   Firstly, we prove that for any
$\xi\in\D(S)$, $\eta=S\xi$ is linearly dependent on $\xi$.
Otherwise, assume that $\eta,\xi$ are linearly independent. Since
$\AA$ is $2-$fold transitive, $\H^{(2)}=[(T\xi,T\eta):\,\,
T\in\AA]=[(T\xi,ST\xi):\,\, T\in\AA]\subseteq \G(S)$. It is a
contradiction. Suppose for $\xi_1,\xi_2$ in $\D(S)$,
$S\xi_1=\lambda_1\xi_1$, $S\xi_2=\lambda_2\xi_2$ and
$S(\xi_1+\xi_2)=\lambda(\xi_1+\xi_2)$. Then
$\lambda_1=\lambda_2=\lambda$. This implies that $S=\lambda I$.\\

``$2\Rightarrow 1$''. Suppose for two linearly independent vectors
$\xi,\eta$ in $\H$, $\{(T\xi,T\eta):\,\, T\in\AA\}\neq \H^{(2)}$. By
lemma 1.1, $\{(T\xi,T\eta):\,\, T\in\AA\}$ is the graph of a closed,
densely defined operator $S$, i.e., $ST\xi=T\eta=TS\xi$. By
assumption, $S=\lambda I$. So $\eta=\lambda \xi$. It is a contradiction.\\

``$1\Rightarrow 3$'' is obvious. ``$3\Rightarrow 1$''. Let $\xi,
\eta$ be linearly independent vectors. Suppose that
$[(T\xi,T\eta):\,\, T\in\AA]\neq \H^{(2)}$, by Lemma 1.1, there is a
closed, densely defined operator $S$ such that
$\G(S)=[(T\xi,T\eta):\,\, T\in\AA]=[(T\xi,ST\xi):\,\, T\in\AA]$. For
any $\zeta\in\H$, by assumption of 3, there exist a sequence $T_n$
in $\AA$ such that $T_n\xi$ converges to $\zeta$ and
$\sup_n\|T_n\eta\|< \infty$. We can assume that $T_n\eta$ weakly
converges to $z$. By Mazur's theorem, there is a sequence $S_n$ such
that $S_n\eta$ is the convex combination of $T_n\eta$ and $S_n\eta$
strongly converges to $z$. Note that $S_n\xi$ strongly converges to
$\zeta$. Since $SS_n\xi=S_n\eta$ and $S$ is a closed, densely
defined operator, $S\zeta=z$. This implies that  $\zeta\in \D(S)$.
Since $\zeta\in\H$ is arbitrary, $\D(S)=\H$. By closed graph
theorem, $S$ is a bounded operator. By Lemma 1.1 and symmetry of
$\xi$ and $\eta$, $S^{-1}$ is a bounded operator on $\H$. Thus $S$
is a bounded, invertible operator with  inverse $S^{-1}$ also a
bounded operator. For any $\lambda\in \mathbb{C}$, consider vectors
$\xi, \lambda\xi+\eta$. Similar arguments show that $S+\lambda I$ is
a bounded, invertible operator with inverse $(S+\lambda I)^{-1}$
also a bounded operator. So $\sigma(S)=\emptyset$. It is a
contradiction. \\

$``4\Rightarrow 1"$. Otherwise, there exist linearly independent
vectors $\xi, \eta$ in $\H$ such that $\G=[(T\xi,T\eta):\,\,
T\in\AA]\neq \H^{(2)}$. Let $P$ be the projection from $\H^{(2)}$
onto subspace $\G$. Then $0< P< I$ and $P\in Lat \AA^{(2)}$. By
assumption, $P\in M_2(\cc)\tensor 1$. Since $ P\neq I$, rank$P=1$.
Therefore, there exists a unitary matrix $U\in M_2(\cc)\tensor \cc
I$ such that $Q=UPU^*=0 \oplus I$.
 Let $\left(\begin{array}{c}
\zeta\\
 \omega
\end{array}\right)=U\left(\begin{array}{c}
\xi\\
 \eta
\end{array}\right)$.  Then
$Q\left(\begin{array}{c}
\zeta\\
 \omega
\end{array}\right)=UPU^*\left(\begin{array}{c}
\zeta\\
 \omega
\end{array}\right)=UP\left(\begin{array}{c}
\xi\\
 \eta
\end{array}\right)=U\left(\begin{array}{c}
\xi\\
 \eta
\end{array}\right)=\left(\begin{array}{c}
\zeta\\
 \omega
\end{array}\right),$
which implies that $\zeta=0$.
 Since  $\xi,\eta$ are linearly independent and $U$ is a unitary operator, $\zeta,
 \omega$ are linearly independent. In particular, $\zeta\neq 0$. It is a
 contradiction.\\

``$1\Rightarrow 4$".  Suppose $\AA$ is $2-$fold transitive. Let
$I_2\neq P=(P_{ij})\in M_{2}(\cc)\tensor\B(\H)$ such that
$P\in\text{Lat}\AA^{(2)}$. It is easy to see
\[P\H^{(2)}=\bigvee\left[\AA^{(2)}\left(\begin{array}{c}
\xi\\
\eta
\end{array}\right): \left(\begin{array}{c}
\xi\\
\eta
\end{array}\right)\in P\H^{(2)}\right].
\]
We only need  to prove that if $P$ is the projection onto space
$\left[\AA^{(2)}\left(\begin{array}{c}
\xi\\
\eta
\end{array}\right)\right]$ for some $\xi, \eta\in\H$,
then $P\in M_2(\cc)\tensor \cc I$.  If $\xi, \eta$ are linearly
independent, by assumption 1, $P=I_2$. If $\xi, \eta$ are linearly
dependent, then $\eta=\lambda\xi$ for some $\lambda\in \mathbb{C}$.
\[\left[\AA^{(2)}\left(\begin{array}{c}
\xi\\
\eta
\end{array}\right)\right]=\left[\left(\begin{array}{c}
T \xi\\
\lambda T\xi\\
\end{array}\right):\,\, T\in\AA\right]=\left[\left(\begin{array}{c}
\zeta\\
\lambda \zeta\\
\end{array}\right):\,\, \zeta\in\H\right]\]
is the closure of the range of the operator
\[\left(\begin{array}{cc} 1&0\\
\lambda&0
\end{array}\right)\in M_2(\cc)\tensor \cc I.\]
Therefore, $P\in M_2(\cc)\tensor \cc I$.

\end{proof}

By using induction on $n$ and similar idea of proof of Proposition
1.2, we can prove the following proposition.
\begin{Proposition}Let $\AA\subseteq\B(\H)$ be a transitive algebra
and $n\in \mathbb{N}$. Then $\AA$ is $n-$fold transitive if and only
if Lat$\AA^{(n)}=$ Lat$\B(\H)^{(n)}$, i.e., for any projection
$P\in\text{Lat}(\cc I_n\tensor\L)$, $P\in M_n(\cc)\tensor \cc I$.
Furthermore, the strong closure of $\AA$ is $\B(\H)$ if and only if
for any $n\in \mathbb{N}$, Lat$\AA^{(n)}=$ Lat$\B(\H)^{(n)}$.
\end{Proposition}

Compare with Arveson's characterizations of $n-$fold transitivity by
graph transformations (see~\cite{Ar}), the characterization of
$n-$fold transitivity given by Proposition 1.3 is more
``computable".

\section{From $n-1$ fold transitivity to $n$ fold transitivity}
Suppose $\AA$ is a transitive algebra on $\H$ and $\AA$ is $n-1$
($n\geq 2$) fold transitive.  By Proposition 1.3, $\AA$ is $n-$fold
transitive if and only if for any projection $P\in\text{Lat}(\cc
I_n\tensor\AA)$, $P\in M_n(\cc)\tensor \cc I$.  It is interesting to
know  under what conditions a projection $P\in\text{Lat}(\cc
I_n\tensor\AA)$ implies that $P\in M_n(\cc)\tensor \cc I$. In the
following, we will provide certain
conditions along this line. \\

 Recall
that for each bounded operator $T\in\B(\H)$, we can associate two
closed subspaces of $\H$, the \emph{null} space $\{\xi\in\H:
T\xi=0\}$ and the \emph{range} space, $[T(\H)]$,  which is the
closure of the range $T(\H)=\{T\xi: \xi\in\H\}$. The corresponding
projections are called the \emph{null projection}, denoted by
$N(T)$, and the \emph{range projection}, $R(T)$, respectively.

\begin{Lemma} Let $P=\left(\begin{array}{cc} T_1 & S\\
S^*& T_2
\end{array}\right)\in \B(\H\bigoplus\K)$ be a projection and
$Z=\left(\begin{array}{cc} X & 0\\
0& Y
\end{array}\right)\in \B(\H\bigoplus\K)$. If $P\in\text{Lat}Z$, then
 $R(T_1)$ and  $N(I-T_1)$ are in $\text{Lat}X$ .
\end{Lemma}
\begin{proof} Since $P$ is a projection, we have
\begin{equation}
T_1(I-T_1)=SS^*,
\end{equation}
which implies that $R(S)\subseteq R(T_1)$. If $P\in\text{Lat}Z$,
then we have equation
\begin{equation}
T_1XT_1+SYS^*=XT_1.
\end{equation}
Therefore, $X(R(T_1))\subseteq R(T_1)$, which implies that
$R(T_1)\in\text{Lat}X$.

$\forall \xi\in N(I-T_1)$, by equation (1), $SS^*\xi=0$, which
implies that $S^*\xi=0$. By equation (2), $T_1X\xi=X\xi$, which
implies that $X\xi\in N(I-T_1)$. So $N(I-T_1)$ is in Lat$X$.

\end{proof}

\begin{Corollary}Let $\AA$ be a transitive algebra. Given a projection
 $P=(P_{ij})_{n\times
n}\in M_{n}(\cc)\tensor\B(\H)$, if $P\in \text{Lat}\AA^{(n)}$, then
$P_{ii}=0$ or $P_{ii}=I$ or $P_{ii}-P_{ii}^2$ is a one-to-one
self-adjoint operator with dense range for $1\leq i\leq n$, i.e,
$N(P_{ii}-P_{ii}^2)=0$ and $R(P_{ii}-P_{ii}^2)=I$.
\end{Corollary}
\begin{proof} For $1\leq i\leq n$,
if $0<P_{ii}<I$, then $R(P_{ii})\neq 0$ and $N(I-P_{ii})\neq I$.
Since $\AA$ is a transitive algebra, by Lemma 2.1, $R(P_{ii})=I$ and
$N(I-P_{ii})=0$. This implies that both $P_{ii}$ and $I-P_{ii}$ are
one-to-one operators. Therefore $P_{ii}(I-P_{ii})$ is a one-to-one
self-adjoint operator. So $N(P_{ii}-P_{ii}^2)=0$ and
$R(P_{ii}-P_{ii}^2)=I$.
\end{proof}

\begin{Lemma} Let $\AA$ be a transitive algebra and $n\geq 2$ be a positive integer.
 Suppose $\AA$ is $(n-1)-$fold
transitive and  $P=(P_{ij})_{n\times n}\in M_{n}(\cc)\tensor\B(\H)$
is a projection in  $Lat\AA^{(n)}$. Write $P
=\left(\begin{array}{cc} T_1 & S\\
S^*& T_2
\end{array}\right)\in\B(\H^{(m)}\bigoplus\H^{(n-m)})$, where $1\leq m\leq n-1$.
 If $R(T_1)\neq
I_{m}$ or $N(I_{m}-T_1)\neq 0_{m}$, then $P\in M_{n}(\cc)\tensor \cc
I$.
\end{Lemma}
\begin{proof}
By Lemma 2.1, $R(T_1), N(I-T_1)\in\text{Lat}\AA^{(m)}$. By
assumption and Proposition 1.3, there exist projections $Q_1, Q_2\in
M_{m}(\cc)\tensor \cc I$ such that $Q_1=R(T_1)$ and
$Q_2=N(I_{m}-T_1)$.\\

 If $R(T_1)\neq I_{m}$, rank$Q_1\leq m-1$.
Since $Q_1\in M_m(\cc)\tensor \cc I$, there exists a unitary
operator $U_1\in M_{m}(\cc)\tensor \cc I$ such that
\[U_1Q_1U_1^*=0_{k}\bigoplus I_{m-k},\]
where $1\leq k\leq m$. Let $W_1= U_1\bigoplus I_{n-m}$. Note that
$W_1(I_{n}\otimes Z)W_1^*=(I_{n}\otimes Z)$, for any $Z\in\AA$.
Therefore, $W_1PW_1^*\in\text{Lat}\AA^{(n)}$. Since $P\in
 M_{n}(\cc)\tensor \cc I$ if and only if $W_1PW_1^*\in
 M_{n}(\cc)\tensor \cc I$, we can assume $R(T_1)=0_{k}\bigoplus
 I_{m-k}.$
This implies that $P_{11}=0$ and therefore $P_{1i}=P_{i1}=0$ for all
$1\leq i\leq n$. So
\[P=\left(\begin{array}{cc} 0 & 0\\
0& P_1
\end{array}\right)\in\text{Lat}\AA^{(n)}\]
and therefore $P_1\in\text{Lat}\AA^{(n-1)}$. By assumption, $P_1\in
M_{n-1}(\cc)\tensor \cc I$ and therefore $P\in M_{n}(\cc)\tensor \cc
I$.\\

If $N(I_{m}-T_1)\neq 0_{m}$, rank$Q_2\geq 1$. Since $Q_2\in
M_m(\cc)\tensor\cc I$, there exists a unitary operator $U_2\in
M_{m}(\cc)\tensor \cc I$ such that
\[U_2Q_2U_2^*=I_{k'}\bigoplus 0_{m-k'},\] where $k'\geq 1$.
Let $W_2=U_2\bigoplus I_{n-m}$. Note that $W_2(I_{n}\otimes
Z)W_2^*=(I_{n}\otimes Z)$, for any $Z\in\AA$. Therefore,
$W_2PW_2^*\in\text{Lat}\AA^{(n)}$. Since $P\in
 M_{n}(\cc)\tensor \cc I$ if and only if $W_2PW_2^*\in
 M_{n}(\cc)\tensor \cc I$, we can assume
 $N(I_{(m)}-T_1)=I_{k'}\bigoplus 0_{m-k'}.$
This implies that $P_{11}=I$  and therefore $P_{1i}=P_{i1}=0$ for
all $1\leq i\leq n$. So
\[P=\left(\begin{array}{cc} I& 0\\
0& P_2
\end{array}\right)\in\text{Lat}\AA^{(n)}\]
and therefore $P_2\in\text{Lat}\AA^{(n-1)}$. By assumption, $P_2\in
M_{n-1}(\cc)\tensor 1$ and therefore $P\in M_{n}(\cc)\tensor \cc I$.
\end{proof}

\begin{Lemma} Let $\AA$ be a transitive algebra and $n\geq 3$ be a positive integer.
 Suppose $\AA$ is $(n-1)-$fold
transitive and  $P=(P_{ij})_{n\times n}\in M_{n}(\cc)\tensor\B(\H)$
is a projection in  $Lat\AA^{(n)}$. Write $P
=\left(\begin{array}{cc} T_1 & S\\
S^*& T_2
\end{array}\right)\in\B(\H\bigoplus\H^{(n-1)})$. If the range
 of $T_1$ is $\H$, then $P\in M_{n}(\cc)\tensor \cc I$.
\end{Lemma}
\begin{proof} Assume $P\notin M_{n}(\cc)\tensor \cc I$. By Corollary
2.2 and Lemma 2.3, both $T_1-T_1^2$ and $T_2-T_2^2$ are one-to-one
operators with dense ranges. By assumption, $T_1$ is an invertible
operator and therefore $T_1^{-1}\in \B(\H)$.  Since $P$ is a
projection, we have $T_1-T_1^2=SS^*$ and $T_2-T_2^2=S^*S$. Let
$S=U|S|$ be the polar decomposition. Then $\H\bigoplus
0^{(n-1)}=R(T_1-T_1^2)=R(S)=UU^*$ and $0\bigoplus
\H^{(n-1)}=R(T_2-T_2^2)=R(S^*)=U^*U$. Hence,
$U=(U_1,\cdots,U_{n-1})$ is a unitary operator  from $0\bigoplus
\H^{(n-1)}$ onto $\H\bigoplus 0^{(n-1)}$ such that
$R(T_1-T_1^2)=UU^*$ and $R(T_2-T_2^2)=U^*U$. Let
$H=\sqrt{T_1-T_1^2}$. It is easy to see
\[P=\left(\begin{array}{cccc}
T_1& HU_1&\cdots&HU_{n-1}\\
U_1^*H&U_1^*(I-T_1)U_1&\cdots&U_1^*(I-T_1)U_{n-1}\\
\vdots&\vdots&\ddots&\vdots\\
U_{n-1}^*H&U_{n-1}^*(I-T_1)U_{2}&\cdots&U_{n-1}^*(I-T_1)U_{n-1}
\end{array}\right)=\]
\[\left(\begin{array}{cccc}
I&&&\\
&U_1^*&&\\
&&\ddots&\\
&&&U_{n-1}^*
\end{array}\right)\left(\begin{array}{cccc}
T_1& H&\cdots&H\\
H&I-T_1&\cdots&I-T_1\\
\vdots&\vdots&\ddots&\vdots\\
H&I-T_1&\cdots&I-T_1
\end{array}\right)\left(\begin{array}{cccc}
I&&&\\
&U_1&&\\
&&\ddots&\\
&&&U_{n-1}
\end{array}\right).\]
For every $0\neq \xi\in \H$,  we have
\[\left(\begin{array}{cc}
T_1&H\\
H&I-T_1
\end{array}\right)\left(\begin{array}{c}
-T_1^{-1}HU_1\xi\\
U_1\xi
\end{array}\right)=\left(\begin{array}{c}
0\\
0
\end{array}\right).\]
Therefore,

\[\left(\begin{array}{cc}
T_1&HU_1\\
U_1^*H&U_1^*(I-T_1)U_1
\end{array}\right)\left(\begin{array}{c}
-U_1^*T_1^{-1}HU_1\xi\\
\xi
\end{array}\right)=\left(\begin{array}{c}
0\\
0
\end{array}\right),\]
which implies that $N\left(\begin{array}{cc}
T_1&HU_1\\
U_1^*H&U_1^*(I-T_1)U_1
\end{array}\right)\neq 0_2$ and therefore $R\left(\begin{array}{cc}
T_1&HU_1\\
U_1^*H&U_1^*(I-T_1)U_1
\end{array}\right)$ $\neq I_2$. By Lemma 2.3, $P\in M_n(\cc)\tensor \cc I$. It is
a contradiction.
\end{proof}

As a corollary of Proposition 1.3, Lemma 2.3 and Lemma 2.4, we get
Foias' result (Proposition 5 of~\cite{Fo}):

\begin{Corollary}\emph{(Foias' Theorem)} Let  $\AA$ be a transitive algebra on a
Hilbert space $\H$. If $\AA$ has no nontrivial invariant operator
ranges, then $\AA$ is strongly dense in $\B(\H)$.
\end{Corollary}

\section{On transitive algebras containing a standard finite von
Neumann subalgebra}

\begin{Theorem} Let $\M$ be a finite von Neumann algebra acting on a
Hilbert space $\H$ and $\AA$ be a transitive algebra which contains
$\M'$. If $\AA$ is $2-$fold transitive, then $\AA$ is strongly dense
in $\B(\H)$.
\end{Theorem}
\begin{proof} We need to prove that $\AA$ is $n-$fold transitive for
$n\in \mathbb{N}$.
 We use induction on $n$. Suppose $\AA$ is $n-$fold transitive
for $n\geq 2$ but $\AA$ is not $n+1$ fold transitive. By Proposition
1.3, let $P=(P_{ij})\in M_{n+1}(\cc)\tensor\B(\H)$ be a projection
such that $P\in\text{Lat}\AA^{(n+1)}$ but $P\notin
M_{n+1}(\cc)\tensor \cc I$. Note that $\M'\subseteq\AA$,
$P\in\text{Lat}\M'^{(n+1)}$. Since $\M'$ is a von Neumann algebra,
$P\in M_{n+1}(\cc)\tensor \M$. Write $P
=\left(\begin{array}{cc} T_1 & S\\
S^*& T_2
\end{array}\right)\in\B(\H\bigoplus\H^{(n)})$. By Corollary 2.2 and Lemma 2.3,
$R(T_1-T_1^2)=\H\bigoplus 0^{(n)}$ and $R(T_2-T_2^2)= 0\bigoplus
\H^{(n)}$. Since $P$ is a projection, we have $T_1(I-T_1)=SS^*$ and
$T_2(I-T_2)=S^*S$. By polar decomposition, there is a unitary
operator $U=(U_2,\cdots,U_{n+1})$ from $0\bigoplus \H^{(n)}$ onto
$\H\bigoplus 0^{(n)}$ such that
\begin{eqnarray}
  UU^* &=& I\bigoplus 0_n, \\
  U^*U &=& 0\bigoplus I_n.
\end{eqnarray}
 By equation (3) and (4),
$I\bigoplus 0_n$ and $0\bigoplus I_n$ are equivalent in
$M_{n+1}(\cc)\tensor \M$. Since $M_{n+1}(\cc)\tensor \M$ is a finite
von Neumann algebra, it is a contradiction.
\end{proof}

\begin{Corollary} Let $\M$ be a standard finite von Neumann algebra on a
Hilbert space $\H$ and $\AA$ be a transitive algebra which contains
$\M$. If $\AA$ is $2-$fold transitive, then $\AA$ is strongly dense
in $\B(\H)$.
\end{Corollary}

\begin{Corollary}\emph{ (Arveson's Theorem)}\,
Let $\AA\subseteq\B(\H)$ be a transitive algebra which contains a
 MASA of $\B(\H)$. Then $\AA$ is strongly dense in  $\B(\H)$.
 \end{Corollary}
\begin{proof} By Theorem 3.1, we need to prove that $\AA$ is
$2-$fold transitive.  Assume $\AA$ is not $2-$fold transitive. By
Proposition 1.2,  there is a non-scalar, closed densely defined
operator $S$ such
  that $TS=ST$ for every $T\in \AA$. Since $\AA$ contains a
  MASA $\A$ of $\B(\H)$, $S$ is affiliated with $\A$. So $S$ is
  a (unbounded)  normal operator. Since $S$ is not a scalar, there is  a non-trivial
  spectral projection $E$ of $S$. By Fuglede's theorem [Fu], $ET=TE$
  for all $T\in\AA$. Since $\AA$ is a transitive algebra, it is a
  contradiction.
\end{proof}

 \begin{Corollary} Let $\M$ be a type $II_1$ factor with a faithful
 normal trace $\tau$. If $\AA$ is a
 transitive subalgebra of $\B(L^2(\M,\tau))$ which contains $\M'$ and a MASA
 of $\M$, then $\AA$ is strongly dense in
 $\B(L^2(\M,\tau))$.
 \end{Corollary}
\begin{proof} Suppose $\A$ is a MASA in $\M$ and $\AA$ is a transitive
algebra which contains $\A$ and $\M'$. By Theorem 3.1, we need to
prove that $\AA$ is $2-$fold transitive. Let $\L$ be the von Neumann
algebra generated by $\A$ and $\M'$. Then $\L'=\A'\cap \M=\A$. Let
$S$ be a closed, densely defined operator such that $ST=TS$ for
every $T\in \AA$. Then $S$ is affiliated with $\A$. Similar
arguments as the proof of Corollary 3.3 show that $S=\lambda I$ for
some $\lambda$. By Proposition 1.2, $\AA$ is $2-$fold transitive.
\end{proof}

\noindent {\bf Remark:} Let $\M$ be a type $II_1$ factor and $\A$ be
a MASA in $\M$. Let $\L$ be the von Neumann subalgebra of
$\B(L^2(\M,\tau))$ generalized by $\A$ and $\M'$. Does $\L$ contains
a MASA of $\B(L^2(\M,\tau))$? It seems that $\M=\L{\mathbb{F}_2}$
and $\A$ be
the MASA generated by one generator is a counterexample. \\

Combine Corollary 3.2 and Proposition 1.2, we have the following
corollary.
\begin{Corollary} Let $\M$ be a type $II_1$ factor with a faithful
normal trace $\tau$. Then the following conditions are equivalent:
\begin{enumerate}
\item for every non scalar, closed, densely defined operator $S$ affiliated with
$\M$, $\{S\}'\cap \B(L^2(\M,\tau))$ has a nontrivial invariant
subspace (affiliated with $\M$);
\item every transitive algebra on $L^2(\M,\tau)$ containing $\M$ is
strongly dense in $\B(L^2(\M,\tau))$.
\end{enumerate}
\end{Corollary}
\medskip

\begin{Lemma} Let $\M$ be a finite von Neumann algebra
and $T\in\M$ be a normal operator, i.e., $TT^*=T^*T$. If $E\in Lat
T\cap\M$, then $ET=TE$. In particular, $E\in Lat T^*$.
\end{Lemma}
\begin{proof} Let $E\in Lat T$. With respect to $E, I-E$, we can
write $T=\left(\begin{array}{cc} A&B\\
0&C \end{array}\right)$. Now we have the following
\[TT^*=\left(\begin{array}{cc} A&B\\
0&C \end{array}\right)\left(\begin{array}{cc} A^*&0\\
B^*&C^* \end{array}\right)=\left(\begin{array}{cc} AA^*+BB^*&BC^*\\
CB^*&CC^* \end{array}\right),\]
\[T^*T=\left(\begin{array}{cc} A^*&0\\
B^*&C^* \end{array}\right)\left(\begin{array}{cc} A&B\\
0&C \end{array}\right)=\left(\begin{array}{cc} A^*A&A^*B\\
B^*A&B^*B+C^*C \end{array}\right).\] By $TT^*=T^*T$,
$AA^*+BB^*=A^*A$. Let $\tau$ be the unique center-valued trace on
$\M$. Apply trace $\tau_E$ induced by $\tau$ on $E\M E$, we get
$\tau_E(BB^*)=0$. This implies that $B=0$ and $ET=TE$.
\end{proof}

\begin{Corollary} Let $\M$ be a finite von Neumann algebra acting on
a Hilbert space $\H$  and $T\in\M$ be a normal operator. If $\AA$ is
an algebra (not necessarily transitive) containing $\M'$ and
$T\in\AA$, then $T^*$ is in the strong closure of $\AA$.
\end{Corollary}
\begin{proof} For $n\in \mathbb{N}$, let $E\in Lat\AA^{(n)}$. Then
$E\in Lat \M'^{(n)}$ and $E\in Lat T^{(n)}$. Hence $E\in Lat
T^{(n)}\cap (M_n(\cc)\tensor \M)$. Since $T^{(n)}$ is a normal
operator in $M_n(\cc)\tensor \M$ and $M_n(\cc)\tensor \M$ is finite,
$E\in Lat T^{*(n)}$ by Lemma 3.6. By Theorem 7.1 of [R-R], $T^*$ is
in the strong closure of $\AA$.
\end{proof}

The following question asked by Radjavi and Rosenthal~\cite{R-R}
remains open: If $\AA$ is a transitive algebra generated by normal
operators, is $\AA$ strongly dense in $\B(\H)$?  But we have the
following result.

\begin{Corollary}Let $\M$ be a type $II_1$ factor with a faithful
 normal trace $\tau$. If $\AA$ is a
 transitive subalgebra of $\B(L^2(\M,\tau))$ generated by normal
 operators in $\M$ and by $\M'$, then  $\AA$ is strongly dense in
 $\B(L^2(\M,\tau))$.
\end{Corollary}
\begin{proof} Suppose $\AA$ is generated by
$\{T_{\alpha}\}\subseteq\M$ and $\M'$, where
$T_{\alpha}T_{\alpha}^*=T_{\alpha}^*T_{\alpha}.$ By Corollary 3.7,
the strong closure of $\AA$ contains the von Neumann algebra $\L$
generated by
 $\{T_{\alpha}\}\subseteq\M$ and $\M'$.  Since $\L$ is  transitive,
 $\L=\B(L^2(\M,\tau))$.
\end{proof}

\section{Non-selfadjoint algebras related to free products of finite von Neumann algebras}

The main results of this section are Theorem 4.1 and Corollary 4.2.
 Examples (Example 4.9) related to free group factors are given at the end of this section.

\begin{Theorem}
Let $(\M,\tau)$ be a finite von Neumann algebra with a faithful
normal trace $\tau$. Suppose $\N\subseteq\M$ is a von Neumann
subalgebra and $Z\in \M$ satisfy the following conditions:
\begin{enumerate}
\item $Z\neq 0$ and $\tau(Z)=0$; \item $Z, Z((\N'\cap\M)\minus
\cc 1), ((\N'\cap\M)\minus \cc 1)Z$ are mutually orthogonal in
$L^2(\M,\tau)$.
\end{enumerate}
Let $\AA\subseteq\B(L^2(\M,\tau))$ be the  algebra generated by $Z,
\N$ and $\M^{'}$. Then $\AA$ is a transitive algebra and strongly
dense in  $\B(L^2(\M,\tau))$.
\end{Theorem}

\begin{Corollary}Let $(\M_1,\tau_1), (\M_2,\tau_2)$ be finite von Neumann algebras and
$\M=(\M_1,\tau_1)*(\M_2,\tau_2)$ be the reduced free product von
Neumann algebra and $\tau$ be the induced faithful normal trace on
$\M$. Suppose $\N\subseteq\M_1$ is a diffuse von Neumann subalgebra
and $Z\in \M_2$ is not a scalar. In $\B(L^2(\M,\tau))$, let $\AA$ be
the algebra generated by $Z$, $\N$ and $\M^{'}$. Then $\AA$ is a
transitive algebra and strongly dense in $\B(L^2(\M,\tau))$.
\end{Corollary}

To prove Theorem 4.1, we need to prove the following auxiliary
lemmas.\\

\begin{Lemma} Let $(\M, \tau)$ be a finite von Neumann algebra
with a faithful normal trace $\tau$ and $\N\subseteq \M$ be a von
Neumann subalgebra. Suppose $Z\in \M$ satisfy $\tau(Z)=0$ and $Z,
Z((\N'\cap\M)\minus \cc 1), ((\N'\cap\M)\minus \cc 1) Z$ are
mutually orthogonal in $L^2(\M, \tau)$, i.e,
$\tau(ZXZ^*)=\tau(Z^*XZ)=\tau(Z^*XZY)=0$ for any $X, Y\in\N'\cap\M$
and $\tau(X)=\tau(Y)=0$. Then $\forall A, B, C, D\in\N'\cap\M$, we
have
\begin{equation}
\tau((AZB)^*(CZD))=\tau(A^*C)\tau(B^*D)\tau(Z^*Z).\end{equation}
\end{Lemma}
\begin{proof} Write $A^*C=X+\tau(A^*C)1$ and $D^*B=Y+\tau(D^*B)1$,
where $X, Y\in \N'\cap\M$ and $\tau(X)=\tau(Y)=0$. Then
\[
\begin{array}{cl}
   & \tau((AZB)^*(CZD)) \\
  = & \tau(B^*Z^*A^*CZD) \\
  = & \tau(Z^*(A^*C)Z(DB^*)) \\
  = & \tau(Z^*(X+\tau(A^*C)1)Z(Y+\tau(DB^*)1)) \\
 = & \tau(Z^*XZY)+\tau(Z^*XZ)\tau(DB^*)+\tau(A^*C)\tau(ZYZ^*)+
\tau(A^*C)\tau(Z^*Z)\tau(DB^*) \\
  = & \tau(A^*C)\tau(Z^*Z)\tau(B^*D). \\
\end{array}\]

\end{proof}

\begin{Corollary} With same assumption of Lemma 4.3, $\forall A, B\in
\N'\cap \M$, $AZB=0$ implies that $A=0$ or $B=0$ or $Z=0$.
\end{Corollary}

\begin{Corollary}With same assumption of Lemma 4.3 and assume
$Z\neq 0$. Then $\forall A, B, C, D\in \N'\cap\M$, $\tau((AZB)^*
(CZD))=0$ if and only if $\tau(A^*C)=0$ or $\tau(B^*D)=0$.
\end{Corollary}

\begin{Lemma} With same assumption of Lemma 4.3 and assume
$Z\neq 0.$ For any positive integer $n$, let $P=(P_{ij})_{2\times
2}\in M_2(\cc)\tensor (\N'\cap\M)$ be a projection. If
$P\in\text{Lat}(I_2\tensor Z)$, then $P\in M_2(\cc)\tensor \cc I$.
\end{Lemma}
\begin{proof}
$P(I_2\tensor Z)P=(I_2\tensor Z)P$ implies the following equation
\begin{equation}\label{1}
P_{11}ZP_{11}+P_{12}ZP_{21}=ZP_{11}.
\end{equation}
For $P_{11}, P_{12}$, by the Gram-Schmidts orthogonal process (with
respect to $\tau$), we have
\begin{eqnarray}
   P_{11}&=&x_1, \\
   P_{12}&=&\alpha_{12}x_1+x_2,
\end{eqnarray}
where $x_1, x_2\in \N'\cap\M$ are  orthogonal in $L^2(\M,\tau)$,
i.e, $\tau(x_2^*x_1)=0$ and $\alpha_{12}\in\cc$. By taking $*$ of
equation (7), (8), we have
\begin{eqnarray}
   P_{11}&=&x_1^*,\\
   P_{21}&=&\bar{\alpha}_{12}x_1^*+x_2^*,
  \end{eqnarray}
where $x_1^*, x_2^*\in \N'\cap\M$ are orthogonal in $L^2(\M,\tau)$.
Plug equations (7),(8),(9),(10) in equation (6), we have
\begin{equation}\label{4}
(1+|\alpha_{12}|^2)x_1Zx_1^*+x_2Zx_2^*+\alpha_{12}x_1Zx_2^*
+\bar{\alpha}_{12}x_2Zx_1^*=Zx_1^*.
\end{equation}

By Corollary 4.5, $x_2Zx_2^*$ is perpendicular with $x_1Zx_1^*$,
$x_1Zx_2^*$,  $x_2Zx_1^*$ and $Zx_1^*$. Therefore, $x_2Zx_2^*=0$. By
Corollary 4.4, we have $x_2=0$. So equation (11) implies that
\[(1+|\alpha_{12}|^2)x_1Zx_1^*=Zx_1^*.\]
By Corollary 4.4, either $x_1=0\in\cc$ or
$x_1=\frac{1}{1+|\alpha_{12}|^2}\in\cc$. Therefore
$P_{11}=x_1\in\cc$, $P_{12}=\alpha_{12}x_1$ and
$P_{21}=P_{12}^*\in\cc$. By symmetry, we have $P_{22}\in\cc$, which
implies that $P\in M_2(\cc)\tensor \cc I$.
\end{proof}

\emph{Proof of Theorem 4.1.}\,\, Let $\L$ be the von Neumann algebra
generated by $\N$ and $\M^{'}$ and  $P\in\B(\H)$ be a projection
such that $P\in\text{Lat}\AA$. Then $P\in \N^{'}\cap \M$. Since
$P\in\text{Lat}Z$, we have $(I-P)ZP=0$. By assumptions of Theorem
4.1 and Corollary 4.4, $P=I$ or $P=0$, which implies that $\AA$ is a
transitive algebra. If $Q\in M_2(\cc)\tensor \B(L^2(\M,\tau))$ is a
projection such that $Q\in\text{Lat}(I_2\tensor\AA)$, then $Q\in
\text{Lat}(I_2\tensor\L)$. Since $\L$ is a von Neumann algebra,
$Q\in (I_2\tensor\L)^{'}=M_2(\cc)\tensor\L^{'}=M_2(\cc)\tensor
(\N^{'}\cap \M).$ Note $Q\in\text{Lat}(I_2\tensor Z)$. By Lemma 4.6,
$Q\in M_2(\cc)\tensor \cc I$. By Proposition 1.2, $\AA$ is $2-$fold
transitive. By Theorem
3.1, $\AA$ is strongly dense in  $\B(L^2(\M,\tau))$. QED\\

 The following lemma is proved by Popa in~\cite{Po}.

\begin{Lemma} Let $\B$ be a diffuse von Neumann
subalgebra of a finite von Neumann algebra $\M$ \emph{(}with a
faithful normal trace $\tau$\emph{)} and $U$ be a unitary operator
in $\M$. If $U\B U^*$ and $\B$ are orthogonal with respect to
$\tau$, then $U$ is orthogonal to $\{V\in\M:\,\, V\B V^*=\B,
VV^*=V^*V=I\}$, the set of normalizers  of $\B$ in $\M$. In
particular, $U$ is orthogonal to $\B'\cap\M$.
\end{Lemma}

The following result is well known. For the sake of completeness, we
include the proof in the following. Our proof follows Popa's idea
in~\cite{Po}.

\begin{Lemma} Let $(\M_1,\tau_1), (\M_2,\tau_2)$ be finite von Neumann algebras and
$\M=(\M_1,\tau_1)*(\M_2,\tau_2)$ be the reduced free product von
Neumann algebra and $\tau$ be the induced faithful normal trace on
$\M$. Suppose $\N$ is a diffuse von Neumann  subalgebra of $\M_1$,
then $\N'\cap \M=\N'\cap \M_1$.
\end{Lemma}
\begin{proof} We only need to prove $\N'\cap \M\subseteq\N'\cap \M_1$.
Let $\dot{\M}_1, \dot{\M}_2$ be the set
$\{T\in\M_1:\,\,\tau_1(T)=0\}$ and  $\{T\in\M_2:\,\,\tau_2(T)=0\}$,
respectively. Note that $L^2(\M,\tau)$ is the closure of $\cc
I\oplus \dot{\M}_1\oplus \dot{\M}_2 \oplus \dot{\M}_1\otimes
\dot{\M}_2\oplus \dot{\M}_2\otimes \dot{\M}_1\oplus \cdots$. To
prove $\N'\cap \M\subseteq\N'\cap \M_1$, we only need to prove
$\dot{\M}_2, \dot{\M}_1\otimes \dot{\M}_2, \dot{\M}_2\otimes
\dot{\M}_1, \cdots$ are orthogonal to $\N'\cap\M_1$. Note that
$\dot{\M}_1$ and $\dot{\M}_2$ are the closure (with respect to
strong operator topology) of linear span of
$\{U\in\M_1:\,\,\tau_1(U)=0, UU^*=U^*U=I\}$ and
$\{V\in\M_2:\,\,\tau_2(V)=0, VV^*=V^*V=I\}$, respectively. We need
to prove the following words: $U_1V_1\cdots U_nV_n$, $U_1V_1\cdots
U_n$ ($n\geq 2$), $V_1U_1V_2U_2\cdots$, are orthogonal to
$\N'\cap\M_1$, where $U_1,U_2,\cdots $ and $V_1,V_2,\cdots$ are
unitary operators in $\dot{\M}_1$ and $\dot{\M}_2$, respectively.
Note that $U_1V_1\cdots U_nV_n \dot{\M}_1 V_n^*U_n^*\cdots V_1^*
U_1^*$ and $\M_1$ are orthogonal in $L^2(\M,\tau)$.  Von Neumann
algebra $U_1V_1\cdots U_nV_n \N V_n^*U_n^*\cdots V_1^* U_1^*$ and
$\N$ are orthogonal with respect to $\tau$. By Lemma 4.7,
$U_1V_1\cdots U_nV_n$ is orthogonal to $\N'\cap\M$. Similarly, we
can prove other cases.
\end{proof}

\emph{Proof of Corollary 4.2.} We can assume $Z\neq 0$ and
$\tau(Z)=0$ (Otherwise, consider $Z-\tau(Z)$).  By Lemma 4.8,
$\N^{'}\cap \M=\N^{'}\cap \M_1$ and therefore $Z$ is free with
$\N^{'}\cap \M$. So $Z, Z((\N'\cap\M)\minus \cc 1),
((\N'\cap\M)\minus \cc 1)Z$ are mutually orthogonal in
$L^2(\M,\tau)$. By Theorem 4.1, $\AA$ is a
transitive algebra and strongly dense in $\B(L^2(\M,\tau))$. QED\\

\noindent{\bf Example 4.9.} \emph{Let $\M=\L{\mathbb{F}_n}$ ($2\leq
n\leq \infty$) be the type $II_1$ factor associated with the left
regular representation $\lambda$ of the free group $\mathbb{F}_n$ on
$n$ generators and $a, b$ be two generators of $\mathbb{F}_n$. Let
$k$ be a positive integer and $Z\in\{\lambda(a),\lambda(a)^*\}''$ be
a nonscalar operator. By Corollary 3.7 and Corollary 4.2, the
algebra $\AA$ on $L^2(\M,\tau)$ generated by $Z, \lambda(b)^k,$ and
$\M'$ is a transitive algebra and strongly dense in
$\B(L^2(\M,\tau))$.}

\section{Non-selfadjoint algebras related to $(M_2(\cc),\frac{1}{2}Tr)*(M_2(\cc),\frac{1}{2}Tr)$}
Let $\M=(M_2(\cc),\frac{1}{2}Tr)*(M_2(\cc),\frac{1}{2}Tr)$ be the
reduced free product von Neumann algebra with respect to the
normalized trace on $M_2(\cc)$. Then $\M$ is a type $II_1$ factor
with a faithful normal trace $\tau$. Let
$\left(E_{ij}\right)_{i,j=1, 2}$ and $\left(F_{ij}\right)_{i,j=1,
2}$ be the matrix units in $(M_2(\cc),\frac{1}{2}Tr)*1$ and
$1*(M_2(\cc),\frac{1}{2}Tr)$, respectively. Consider the following
subalgebras of $\B(L^2(\M,\tau))$:
\begin{enumerate}
\item $\AA_1$, the algebra generated by $E_{11}, E_{12}, F_{11}$ and
$\M'$;
\item $\AA_2$, the algebra generated by $E_{11}, E_{12}, F_{12}$ and
$\M'$;
\item $\AA_3$, the algebra generated by $E_{11}, F_{12}$ and $\M'$;
\item $\AA_4$, the algebra generated by $E_{12}, F_{12}$ and $\M'$.
\end{enumerate}

We have the following questions: For $1\leq i\leq 4$, is $\AA_i$ a
transitive algebra? If $\AA_i$ is a transitive algebra, is $\AA_i$
strongly dense in $\B(L^2(\M,\tau)))$?  In this section we  will
prove the following results.
\begin{enumerate}
\item $\AA_1$ is a transitive algebra and  strongly dense in
$\B(L^2(\M,\tau))$.
\item $\AA_3,\AA_4$ are not transitive. The invariant subspaces of
$\AA_3$ and $\AA_4$ are abundant.
\end{enumerate}
The question on $\AA_2$ remains open!\\

To prove our main results, we need to introduce the following
operators. In $(M_2(\cc),\frac{1}{2}Tr)*1$, let
\[W_0=\left(\begin{array}{cc}
1&0\\
0&1
\end{array}\right), W_1=\left(\begin{array}{cc}
1&0\\
0&-1
\end{array}\right), W_2=\left(\begin{array}{cc}
0&-1\\
1&0
\end{array}\right), W_3=\left(\begin{array}{cc}
0&1\\
1&0
\end{array}\right).\]
Then we have
\begin{enumerate}
\item $W_1^2=W_2^2=W_3^2=W_0$; \item $W_1W_2=-W_2W_1=-W_3$,
$W_1W_3=-W_3W_1=-W_2$, $W_2W_3=-W_3W_2=-W_1$; \item $W_0, W_1, W_2,
W_3$ form an orthonormal base of $L^2((M_2(\cc),\frac{1}{2}Tr)*1,
\tau)$.
\end{enumerate}
Similarly, in $1*(M_2(\cc),\frac{1}{2}Tr)$, let
\[V_0=\left(\begin{array}{cc}
1&0\\
0&1
\end{array}\right), V_1=\left(\begin{array}{cc}
1&0\\
0&-1
\end{array}\right), V_2=\left(\begin{array}{cc}
0&-1\\
1&0
\end{array}\right), V_3=\left(\begin{array}{cc}
0&1\\
1&0
\end{array}\right).\]
Then we have similar equations and $\M$ is the von Neumann algebra
generated by $\{W_0,W_1,\\ W_2,W_3\}* \{V_0,V_1, V_2,V_3\}$.

\subsection{The case of $\AA_1$}

 Let
$W=W_1$, $V=V_1$ and $U=WV$. Then $U$ is a haar unitary operator in
$\M$. Let $\A$ be the von Neumann algebra generated by $U$ and
$U^*$. The following lemma is an easy exercise.

\begin{Lemma} $\A$ is a MASA of $\M$.
\end{Lemma}

\begin{Theorem} Let $(\M,\tau)$ and $\A$ be as above. For $Z\in
 (M_2(\cc),\frac{1}{2}Tr)*1$, let $\AA$ be
the algebra generated by $Z, \A$ and $\M^{'}$ in $\B(L^2(\M,\tau))$.
Then $\AA$ is a transitive algebra if and only if $Z\notin\{W\}''$.
In this case, $\AA$ is strongly dense in $\B(L^2(\M,\tau))$.
\end{Theorem}

\begin{Corollary} $\AA_1$ is a transitive algebra and strongly dense
in $\B(L^2(\M,\tau))$.
\end{Corollary}

We need the following lemmas.
\begin{Lemma} With the same assumption of Theorem 5.2 and
assume $Z=\lambda W$ for some $\lambda\in\cc$. Then $\AA$ is not
transitive.
\end{Lemma}
\begin{proof}
We can assume $Z=W$. Then $\AA$ is a von Neumann algebra. Note
\[(U+U^*)W=(WV+VW)W=WVW+V=WVW+W^2V=W(U+U^*).\]
So $U+U^*\in\AA^{'}$, which implies that $\AA$ is not transitive.
\end{proof}

\begin{Lemma} With the same assumption of Theorem 5.2 and
assume $Z\neq 0$ and $\tau(Z)=\tau(Z^*W)=0$. Then $\AA$ is a
transitive algebra and strongly dense in $\B(L^2(\M,\tau))$.
\end{Lemma}
\begin{proof} Note $\tau(Z^*V)=0$ and $W^*=W$. It is easy to verify
\[\tau(Z(WV)^nZ^*)=\tau(Z^*(WV)^nZ)=\tau((WV)^nZ(WV)^mZ^*)=0\] for
all $m, n\in \zz\diagup \{0\}$. Since $\A\ominus \mathbb{C}I$ is
generated by $(WV)^n$ for all $n\in\zz$ and $n\neq 0$, $Z,
Z(\A\minus \cc 1), (\A\minus \cc 1)Z$ are mutually orthogonal in
$L^2(\M,\tau)$.
 Since $\A$ is a MASA of $\M$,
$\A^{'}\cap \M=\A$. By Theorem 4.1, $\AA$ is a transitive algebra
and strongly dense in  $\B(L^2(\M,\tau))$.
\end{proof}

\begin{Lemma} With the same assumption of Theorem 5.2 and
assume $Z=W+X$, $X\neq 0$ and $\tau(X)=\tau(X^*W)=0$. Then $\AA$ is
transitive.
\end{Lemma}
\begin{proof}
Let $P\in \B(L^2(\M,\tau))$ be a projection such that $P\in Lat\AA$.
Since $\AA$ contains $\A$ and $\M'$,  $P\in
\{\A,\M'\}'=\A'\cap\M=\A$. By identifying $\M$ as a subset of
$L^2(\M,\tau)$, we can write $P=\sum_{n=-\infty}^{\infty}\lambda_{n}
U^n$, where  $\lambda_{n}=\overline{\lambda}_{-n}.$ Since $P\in Lat
Z$, we have $PZP=ZP.$ Thus, we have the following equation.

\[\sum_{m,n\in\zz} \lambda_{m}\lambda_{n}U^mWU^n
+\sum_{m,n\in\zz} \lambda_{m}\lambda_{n}U^mXU^n= \sum_{n\in\zz}
\lambda_{n}WU^n+ \sum_{n\in\zz}\lambda_{n}XU^n.\]

It is easy to verify that $U^{-n}XU^n$, $n\in\zz\diagup\{0\}$ are
mutually orthogonal in $L^2(\M,\tau)$ and perpendicular to $U^lWU^m,
U^rXU^s, WU^m, XU^m$ for all $l,m, r,s\in\zz$ and $r\neq -s$.
Therefore, for $n\neq 0$, the coefficients of $U^{-n}XU^n$ are zero,
i.e, $|\lambda_n|^2=0$, which implies that $\lambda_n=0, \forall
n\in\zz\diagup\{0\}$. So $P=\lambda_0 I$. Since $P=P^*=P^2$, $P=0$
or $I$. This implies that $\AA$ is transitive.
\end{proof}

\emph{Proof of Theorem 5.2.}\,\, If $Z=\lambda I$, then $\AA$ is the
von Neumann algebra generated by $\A$ and $\M'$. So $\AA$ is not
transitive. If $Z\in\{W\}''$ and $Z\neq \lambda I$, then $\AA$ is
the algebra generated by $W, \AA$ and $\M'$.  By Lemma 5.4, $\AA$ is
not transitive. Assume that $Z\notin\{W\}''$. We can assume that
$Z\neq 0$ and $\tau(Z)=0$. Write $Z=\alpha W+X$, where $X\neq 0$ and
$\tau(X)=\tau(X^*W)=0$. If $\alpha=0$, by Lemma 5.5, $\AA$ is a
transitive algebra and  strongly dense in $\B(L^2(\M,\tau))$. If
$\alpha\neq 0$, we can assume $\alpha=1$. By Lemma 5.6, $\AA$ is
transitive. By Lemma 5.1 and Corollary 3.4, $\AA$ is strongly dense
in  $\B(L^2(\M,\tau))$. QED

\subsection{Brown measures of $W_1+F_{12}$ and $E_{12}+F_{12}$}
The following observations are crucial to prove that $\AA_3$ and
$\AA_4$ are not transitive: $\AA_3\subseteq \{(W_1+F_{12})^2\}'$ and
$\AA_4\subseteq \{(E_{12}+F_{12})^2\}'$. So we only need to prove
that the support of Brown measures of $(W_1+F_{12})^2$ and
$(E_{12}+F_{12})^2$ are not single point and then conclude the
existence of nontrivial hyperinvariant subspaces of $(W_1+F_{12})^2$
and $(E_{12}+F_{12})^2$ by~\cite{H-S}.
\subsubsection{Brown measures of $R-$diagonal operators}

Let $\M$ be a finite von Neumann algebra with a faithful normal
trace $\tau$. The Fuglede-Kadison determinant (~\cite{F-K}) of
$T\in\M$ is defined by
\[\Delta(T)=\exp[\tau(\ln(T^*T)^{\frac{1}{2}})]\]
is a generalization of a determinant of a finite matrix.\\

The Brown measure (~\cite{Br}) $\mu_T$ of the element $T$ is a
Schwartz distribution on the complex plane defined by
\[\mu_T=\frac{1}{2\pi}(\frac{\partial^2}{\partial
x^2}+\frac{\partial^2}{\partial y^2})\ln\Delta[T-(x+iy)I].\] If
$\M=M_n(\cc)$ and $\tau=\frac{1}{n} Tr$ is the normalized trace on
$M_n(\cc)$, then $\mu_T$ is the normalized counting measure
$\frac{1}{n}\left(\delta_{\lambda_1}+\delta_{\lambda_2}\cdots+\delta_{\lambda_n}\right)$,
 where $\lambda_1,\lambda_2 \cdots,\lambda_n$ are the eigenvalues of $T$
 repeated according to root multiplicity. If $T$ is normal, $\mu_T$
 is the trace $\tau$ composed with the spectral measure for $T$.
 From the definition, Brown measure $\mu_T$ only depends on the joint
 distribution of $T$ and $T^*$.\\

The Brown measure has the following properties (see~\cite{Br}):
$\mu_T$ is the unique compactly supported measure on $\cc$ such that
$\ln\Delta[T-(x+iy)I]=\int_{\cc} \ln|z-\lambda|d\mu_T(z)$ for all
$\lambda\in\cc$. The support of $\mu_T$ is contained in $\sigma(T)$,
the spectra of $T$. $\mu_{AB}=\mu_{BA}$ for arbitrary $A,B$ in $\M$,
and if $f(z)$ is analytic in a neighborhood of $\sigma(A)$,
$\mu_{f(T)}=(\mu_T)_f$. If $E\in\M$ is a projection such that $E\in
Lat T$, then with respect to $E, I-E$ we can write
\[T=\left(\begin{array}{cc}
A&B\\
0&C
\end{array}\right),\]
where $A=ETE$ and $C=(I-E)T(I-E)$ are elements of the $II_1$ factor
$\M_1=E\M E$ and $\M_2=(I-E)\M (I-E)$. Let $\mu_{A}$ and $\mu_{C}$
be the Brown measure of $A$ and $C$ computed relative to $\M_1$ and
$\M_2$, respectively, then $\mu_T=\alpha\mu_A+(1-\alpha)\mu_C$,
where $\alpha=\tau(E)$. \\

$R-$diagonal operators are introduced by Nica and Speicher
in~\cite{N-S}. Recall that an operator $T$ in a non-commutative
probability space is an $R-$diagonal operator if the $R-$transform
$R_{\mu(T,T^*)}$ of the joint distribution $\mu(T,T^*)$ of $T, T^*$
is of the form
\[R_{\mu(T,T^*)}(z_1,z_2)=\sum_{n=1}^{\infty}\alpha_n (z_1z_2)^n+
\sum_{n=1}^{\infty}\alpha_n (z_2z_1)^n.\] Nica and
Speicher~\cite{N-S} proved that $T$ is an $R-$diagonal operator if
and only if $T$ has same *-distribution as product $UH$, where $U$
and $H$ are *-free random variables in some tracial non commutative
C*-probability space, $U$ is a Haar unitary operator and $H$ is
positive. In~\cite{H-L}, Haagerup and Larson explicitly computed the
Brown measures of \emph{$R-$diagonal} operators in a finite von
Neumann algebra.
\begin{Theorem} \emph{(Theorem 4.4 of~\cite{H-S})} Let $U,H$ be *-free random variables in $(\M,\tau)$,
with $U$ a Haar unitary operator and $H$ a positive operator such
that the distribution $\mu_H$ of $H$ is not a Dirac measure. Then
the Brown measure $\mu_{UH}$ of $UH$ can be computed as the
following.
\begin{enumerate}
\item $\mu_{UH}$ is rotation invariant and its support is the
annulus with inner radius $\|H^{-1}\|_2^{-1}$ and outer radius
$\|H\|_2$.
\item $\mu_{UH}(\{0\})=\mu_{H}(\{0\})$ and for $t\in ]\mu_H(\{0\}),
1]$,
\[\mu_{UH}\left(\mathbf{B}\left(0, \left(\S_{\mu_{H^2}}(t-1)\right)^{-1/2}\right)\right)=t,\]
where $\S_{\mu_{H^2}}$ is the $\S-$transform of $H^2$ and
$\mathbf{B}(0,r)$ is the closed disc with center 0 and radius $r$;
\item $\mu_{UH}$ is the only rotation invariant symmetric
probability measure satisfying 2.
\end{enumerate}
\end{Theorem}

Based on Theorem 5.7, the Brown measure of the sum of a random
variable with an arbitrary distribution and a free $R-$diagonal
element, e.g, $U_n+U_{\infty}$, where $U_n$ and $U_{\infty}$ are the
generators of $\mathbb{Z}_n$ and $\mathbb{Z}$ respectively in the
free product $\mathbb{Z}_n*\mathbb{Z}$, is computed by Biane and
Lehner in~\cite{B-L}.

\subsubsection{Brown measure of $W_1F_{12}$}
\begin{Lemma} $T=W_1F_{12}$ is an $R-$diagonal operator with Brown
measure $\mu_T$ satisfying:
\begin{enumerate}
\item $\mu_T$ is rotation invariant and
 the support of $\mu_T$ is $\mathbf{B}\left(0,\frac{1}{\sqrt{2}}\right)$;
\item $d\mu_T(z)=\frac{1}{2}\delta_{0}+\frac{1}{2\pi}\frac{1}{(1-r^2)^2}drd\theta$
for $z=re^{i\theta}$ and $0<r\leq \frac{1}{\sqrt{2}}$.
\end{enumerate}
\end{Lemma}
\begin{proof} Note that $F_{12}=V_2^*F_{22}$. Since $W_1$ and $V_2$
are free unitary operator such that $\tau(W_1)=\tau(V_2)=0$, simple
computations show that $\widetilde{U}=W_1V_2^*$ is a Haar unitary
operator. To prove $T$ is an $R-$diagonal operator, we only need to
check that $\widetilde{U}$ is
*-free with $F_{22}$. Note that $F_{22}$ is in the algebra generated
by $V_1$, we only need to prove that $\widetilde{U}$ is *-free with
$V_1$. Note that $V_1^2=I$ and $\tau(V_1)=0$. We need to check
$\tau(\widetilde{U}^{n_1}V_1\widetilde{U}^{n_2}V_1\widetilde{U}^{n_3}\cdots
V_1\widetilde{U}^{n_k})=0$, where $n_1,n_2,\cdots, n_k\in
\mathbb{Z}$ and $n_2,\cdots,n_{k-1}$ are not equal 0. Consider the
word $\widetilde{U}^l V_1\widetilde{U}^m$, $l, m\neq 0$. We have the
following four cases:
\begin{enumerate}
\item If $l,m > 0$, then $\widetilde{U}^l V_1\widetilde{U}^m= W_1V_2^*\cdots
W_1(V_2^*V_1)W_1V_2^*\cdots W_1V_2^*;$
\item If $l>0, m<0$, then $\widetilde{U}^l V_1\widetilde{U}^m= W_1V_2^*\cdots
W_1(V_2^*V_1V_2^*)W_1V_2^*\cdots W_1;$
\item If $l<0, m>0$, then $\widetilde{U}^l V_1\widetilde{U}^m= V_2^*W_1\cdots
W_1V_1W_1V_2^*\cdots W_1V_2^*;$
\item If $l<0, m<0$, then $\widetilde{U}^l V_1\widetilde{U}^m= V_2^*W_1\cdots
W_1(V_1V_2^*)W_1V_2^*\cdots W_1.$
\end{enumerate}
 Note that
$\tau(V_2^*V_1)=\tau(V_2^*V_1V_2^*)=\tau(V_1V_2^*)=0$.
$\widetilde{U}^{n_1}V_1\widetilde{U}^{n_2}V_1\widetilde{U}^{n_3}\cdots
V_1\widetilde{U}^{n_k}$ is an alternating product of centered
elements from $(M_2(\cc),\frac{1}{2}Tr)*1$ and
$1*(M_2(\cc),\frac{1}{2}Tr)$. Thus
$\tau(\widetilde{U}^{n_1}V_1\widetilde{U}^{n_2}V_1\widetilde{U}^{n_3}\cdots
V_1\widetilde{U}^{n_k})=0$. This proves that $T$ is an $R-$diagonal
operator.\\

Now we apply Theorem 5.7 to compute the Brown measure of $T$.  Note
that $H=F_{22}$ and
$\|F_{22}\|_2=(\tau(F_{22}^2))^{1/2}=\frac{1}{\sqrt{2}}$.  Since the
kernel of $F_{22}$ is nontrivial, $\|F_{22}^{-1}\|_2^{-1}=0$. By 1
of Theorem 5.7, we obtain 1 of Lemma 5.8. By 2 of Theorem 5.7,
$\mu_{T}(0)=\mu_H(0)=\frac{1}{2}$.  To compute the density function
of Brown measure of $T$, we first compute the $\S-$transform of
$H^2=H=F_{22}$. Simple computation shows that
$\S_{\mu_{H^2}}(\omega)=\frac{2(\omega+1)}{2\omega+1}$. By 2 of
Theorem 5.7, $t=\mu_{T}\left(\mathbf{B}\left(0,
\left(\S_{\mu_{H^2}}(t-1)\right)^{-1/2}\right)\right)=\mu_{T}\left(\mathbf{B}\left(0,
\sqrt{1-\frac{1}{2t}}\right)\right)$ for $t\in ]\frac{1}{2}, 1]$.
Let $r=\sqrt{1-\frac{1}{2t}}$. Then $t=\frac{1}{2(1-r^2)}$ and
$\mu_{T}\left(\mathbf{B}\left(0, r\right)\right)=\frac{1}{2(1-r^2)}$
for $0<r<\frac{1}{\sqrt{2}}$. This implies 2 of Lemma 5.8.
\end{proof}

The following lemma is useful to compute the Brown measures of
$E_{12}+F_{12}$ and $W_1+F_{12}$.

\begin{Lemma} Let $\N=1*(M_2(\cc),\frac{1}{2}Tr)$. Then $W_1\N W_1$ is free with
$\N$.
\end{Lemma}
\begin{proof} Consider an alternating product of  elements of
$W_1\N W_1$ and $\N$: $A_1(W_1B_1W_1)A_2(W_2B_2\\ W_2)A_3 \cdots
(W_1B_{n-1}W_1)A_n$, where $A_2,\cdots,A_{n-1}, B_1,\cdots,B_{n-1}$
are centered elements in $\N$ and $A_1, A_n$ are either centered
elements or 1. Then it an alternating product of elements of
$(M_2(\cc),\frac{1}{2}Tr)*1$ and $\N=1*(M_2(\cc),\frac{1}{2}Tr)$.
Thus the trace is 0.
\end{proof}

Since the Brown measure of $T$ only depends on the joint
distribution of $T$ and $T^*$, the Brown measures of $E_{12}+F_{12}$
 and $W_1+F_{12}$ are same as the Brown measures of
$W_1F_{12}W_1+F_{12}$ and $W_1V_1W_1+F_{12}$, respectively.\\

With respect to matrix units of $\N$, write
$W_1=\left(\begin{array}{cc} A&B^*\\
B&C\end{array}\right)$. Since $W_1^2=I$, we have
\[\left(\begin{array}{cc} 1&0\\
0&1\end{array}\right)=\left(\begin{array}{cc} A&B^*\\
B&C\end{array}\right)\left(\begin{array}{cc} A&B^*\\
B&C\end{array}\right)=\left(\begin{array}{cc}
A^2+B^*B& AB^*+B^*C\\
BA+CB&C^2+BB^*\end{array}\right).\] Hence,
\begin{equation}
-CB=BA.
\end{equation}

\subsubsection{Brown measure of  $E_{12}+F_{12}$}
 Since $\mu_{(E_{12}+F_{12})}=\mu_{(W_1F_{12}W_1+F_{12})}$, we
  need to compute the Brown measure of $W_1F_{12}W_1+F_{12}$. Note
that
$(W_1F_{12}W_1+F_{12})^2=W_1F_{12}W_1F_{12}+F_{12}W_1F_{12}W_1= \left(\begin{array}{cc} A&B^*\\
B&C\end{array}\right)\left(\begin{array}{cc} 0&I\\
0&0\end{array}\right)\left(\begin{array}{cc} A&B^*\\
B&C\end{array}\right)\\
\left(\begin{array}{cc} 0&I\\
0&0\end{array}\right)+
\left(\begin{array}{cc} 0&I\\
0&0\end{array}\right)\left(\begin{array}{cc} A&B^*\\
B&C\end{array}\right)\left(\begin{array}{cc} 0&I\\
0&0\end{array}\right)\left(\begin{array}{cc} A&B^*\\
B&C\end{array}\right)=\left(\begin{array}{cc} B^2&AB+BC\\
0&B^2\end{array}\right).$ So the Brown measure of
$(W_1F_{12}W_1+F_{12})^2$ is same as the Brown measure of $B^2$ and
the Brown measure of $W_1F_{12}W_1+F_{12}$ is same as the Brown
measure of $B$ (because Brown measures of $W_1F_{12}W_1+F_{12}$ and
$B$ are both rotation invariant). Now we only need to compute the
Brown measure of $B$. Note that
\[W_1F_{12}=\left(\begin{array}{cc} 0&A\\
0&B\end{array}\right).\]  Thus
$\mu_{W_1F_{12}}=\frac{1}{2}\delta_0+\frac{1}{2}\mu_B$. By Lemma
5.7, we conclude that  $\mu_B$ is rotation invariant and
 the support of $\mu_B$ is $\mathbf{B}\left(0,\frac{1}{\sqrt{2}}\right)$
 and $d\mu_B(z)=\frac{1}{\pi}\frac{1}{(1-r^2)^2}drd\theta$ for
$z=re^{i\theta}$ and $0\leq r\leq \frac{1}{\sqrt{2}}$. Summarize
above, we have the following proposition.
\begin{Proposition} Let $\mu_{(E_{12}+F_{12})}$ be the Brown measure
of $E_{12}+F_{12}$. Then we have the following
\begin{enumerate}
\item $\mu_{(E_{12}+F_{12})}$ is rotation invariant and
 the support of $\mu_{(E_{12}+F_{12})}$ is
 $\mathbf{B}\left(0,\frac{1}{\sqrt{2}}\right)$;
\item $d\mu_{(E_{12}+F_{12})}(z)=\frac{1}{\pi}\frac{1}{(1-r^2)^2}drd\theta$ for
$z=re^{i\theta}$ and $0\leq r\leq \frac{1}{\sqrt{2}}$.
\end{enumerate}
\end{Proposition}
\begin{Corollary}Let $\mu_{(E_{12}+F_{12})^2}$ be the Brown measure
of $(E_{12}+F_{12})^2$. Then we have the following
\begin{enumerate}
\item $\mu_{(E_{12}+F_{12})^2}$ is rotation invariant and
 the support of $\mu_{(E_{12}+F_{12})^2}$ is
 $\mathbf{B}\left(0,\frac{1}{2}\right)$;
\item $d\mu_{(E_{12}+F_{12})^2}(z)=\frac{1}{4\pi}\frac{1}{r(1-r)^2}drd\theta$ for
$z=re^{i\theta}$ and $0\leq r\leq \frac{1}{2}$.
\end{enumerate}
\end{Corollary}

\subsubsection{Brown measure of $W_1+F_{12}$}
 Since $\mu_{(W_{1}+F_{12})}=\mu_{(W_1V_1W_1+F_{12})}$, we
  need to compute the Brown measure of $W_1V_1W_1+F_{12}$.
 We compute the Brown measures of
$(W_1V_1W_1+F_{12})^2$ first. Note that
$(W_1V_1W_1+F_{12})^2=I_2+W_1V_1W_1F_{12}+F_{12}W_1V_1W_1=I_2+ \left(\begin{array}{cc} A&B^*\\
B&C\end{array}\right)\left(\begin{array}{cc} I&0\\
0&-I\end{array}\right)\left(\begin{array}{cc} A&B^*\\
B&C\end{array}\right)\left(\begin{array}{cc} 0&I\\
0&0\end{array}\right)+
\left(\begin{array}{cc} 0&I\\
0&0\end{array}\right)
\left(\begin{array}{cc} A&B^*\\
B&C\end{array}\right)\left(\begin{array}{cc} I&0\\
0&-I\end{array}\right)\left(\begin{array}{cc} A&B^*\\
B&C\end{array}\right)=\left(\begin{array}{cc} I+2BA&A^2-B^*B+BB^*-C^2\\
0&I+2BA\end{array}\right),$ the last equation follows by (12). So
the Brown measure of $(W_1F_{12}W_1+F_{12})^2$ is same as the Brown
measure of $1+2BA$. Now we  only need to compute the Brown measure
of $2BA$. Note that
\[W_1V_1W_1F_{12}=\left(\begin{array}{cc}0&A^2-B^*B\\
0&2BA\end{array}\right).\]   Since $W_1V_1W_1$ is
*-free with $F_{12}$, $\mu_{W_1F_{12}}=\mu_{W_1V_1W_1F_{12}}=
\frac{1}{2}\delta_0+\frac{1}{2}\mu_{2BA}$. By Lemma 5.7, we conclude
that  $\mu_{2BA}$ is rotation invariant and
 the support of $\mu_{2BA}$ is $\mathbf{B}\left(0,\frac{1}{\sqrt{2}}\right)$
 and $d\mu_{2BA}(z)=\frac{1}{\pi}\frac{1}{(1-r^2)^2}drd\theta$ for
$z=re^{i\theta}$ and $0\leq r\leq \frac{1}{\sqrt{2}}$. Summarize
above, we have the following proposition.
\begin{Proposition} Let $\mu_{(W_1+F_{12})^2}$ be the Brown measure
of $(W_1+F_{12})^2$. Then we have the following
\begin{enumerate}
\item $\mu_{(W_1+F_{12})^2}$ is rotation invariant with respect to 1 and
 the support of $\mu_{(W_1+F_{12})^2}$ is
 $\mathbf{B}\left(1,\frac{1}{\sqrt{2}}\right)$;
\item $d\mu_{(W_1+F_{12})^2}(z)=\frac{1}{\pi}\frac{1}{(1-r^2)^2}drd\theta$ for
$z-1=re^{i\theta}$ and $0\leq r\leq \frac{1}{\sqrt{2}}$.
\end{enumerate}
\end{Proposition}

Since the joint $*-$distribution of $W_1+F_{12}, (W_1+F_{12})^*$ and
the joint $*-$distribution of $-(W_1+F_{12}), -(W_1+F_{12})^*$ are
same. We have the following corollary.

\begin{Corollary}Let $\mu_{(W_1+F_{12})}$ be the Brown measure
of $W_1+F_{12}$. Then
 the support of $\mu_{(W_1+F_{12})}$ is $\left\{z\in \mathbb{C}:\,\,
 |z^2-1|\leq \frac{1}{\sqrt{2}}\right\}$.
 \end{Corollary}

\subsection{The case of $\AA_3$}
\begin{Theorem} $\AA_3$ is not a transitive algebra. Indeed, for any
$0\leq r\leq 1$, there is a projection $E_r\in\M$ such that
$\tau(E_r)=r$ and $E_r\in Lat \AA_3$.
\end{Theorem}
\begin{proof} $\AA_3$ is the algebra generated by $W_1, F_{12}$ and
$\M'$. Simple computation shows that
$W_1(W_1+F_{12})^2=(W_1+F_{12})^2 V_1$ and
$F_{12}(W_1+F_{12})^2=(W_1+F_{12})^2F_{12}$. So $\AA_3\subseteq
\{(W_1+F_{12})^2\}'$. Let $0\leq r\leq 1$. By Proposition 5.12,
there is $s$, $0\leq s\leq \frac{1}{\sqrt{2}}$ such that
$\mu_{(W_1+F_{12})^2}\left(\mathbf{B}\left(1,s\right) \right)=r$. By
Theorem 7.1 of~\cite{H-S}, there is a hyperinvariant subspace $E_r$
in $\M$ of $(W_1+F_{12})^2$ such that $\tau(E_r)=r$.
\end{proof}
\subsection{The Case of $\AA_4$}
\begin{Theorem} $\AA_4$ is not a transitive algebra. Indeed, for any
$0\leq r\leq 1$, there is a projection $F_r\in\M$ such that
$\tau(F_r)=r$ and $F_r\in Lat \AA_4$.
\end{Theorem}
\begin{proof} $\AA_4$ is the algebra generated by $E_{12}, F_{12}$ and
$\M'$. Simple computation shows that
$E_{12}(E_{12}+F_{12})^2=(E_{12}+F_{12})^2 E_{12}$ and
$F_{12}(E_{12}+F_{12})^2=(E_{12}+F_{12})^2F_{12}$. So
$\AA_4\subseteq \{(E_{12}+F_{12})^2\}'$. By Corollary 5.11, there is
$s$, $0\leq s\leq \frac{1}{2}$ such that
$\mu_{(E_{12}+F_{12})^2}\left(\mathbf{B}\left(1,s\right) \right)=r$.
By Theorem 7.1 of~\cite{H-S}, there is a
 hyperinvariant subspace $F_r$ in
$\M$ of $(E_{12}+F_{12})^2$ such that $\tau(F_r)=r$.
\end{proof}

\subsection{On the case of $\AA_2$}
Let $\widetilde{\AA}_2\subseteq \B(L^2(\M,\tau))$ be the algebra
generated by $(M_2(\cc),\frac{1}{2}Tr)*1$, $F_{12}$ and $\M'$. Then
$\AA_2\subseteq \widetilde{\AA}_2$. By the following proposition, we
 only need to consider if $\widetilde{\AA}_2$ is transitive and the strong closure
of $\widetilde{\AA}_2$ is $\B(L^2(\M,\tau))$ or not.
\begin{Proposition} $\AA_2$ is transitive if and only if $\widetilde{\AA}_2$ is
transitive; the strong closure of $\AA_2$ is $\B(L^2(\M,\tau))$ if
and only if the strong closure of $\widetilde{\AA}_2$ is
$\B(L^2(\M,\tau))$.
\end{Proposition}
\begin{proof} Suppose $\widetilde{\AA}_2$ is transitive. Let $P\in Lat\AA_2$.
Then $P\in\M$. With respect to matrix units of
$(M_2(\cc),\frac{1}{2}Tr)*1$, $\M\cong M_2(\mathbb{C})\tensor \M_1$.
In the following, we identify $\M$ with $M_2(\cc)\tensor \M_1$.
 Since $PE_{11}=E_{11}P$, $P=\left(\begin{array}{cc}
P_1&0\\
0& P_2
\end{array}\right)$, where $P_1,P_2\in\M_1$.  By $PE_{12}P=E_{12}P$ and simple computation,
we have $P_2\leq P_1$.
 Write $F_{12}=\left(\begin{array}{cc} A&B\\
C&D
\end{array}\right)$, where $A,B,C,D\in\M_1$. By $PF_{12}P=F_{12}P$ and simple computation, we have
$(I-P_2)CP_1=0$. It can be proved that $C$ is a one-to-one operator
with dense range. So $\tau_{\M_1}(P_1)=\tau_{\M_1}(R(CP_1))\leq
\tau_{\M_1}(P_2)$, where $\tau_{\M_1}$ is the trace induced by
$\tau$ on $\M_1$. Since $P_2\leq P_1$, this implies that $P_1=P_2$.
Therefore, $P\in Lat\widetilde{\AA}_2$. Since $\widetilde{\AA}_2$ is
transitive, $P=0$ or
$I$. So $\AA_2$ is transitive.\\

Now suppose the strong closure of $\widetilde{\AA}_2$ is
$\B(L^2(\M,\tau))$. Then $\widetilde{\AA}_2$ is transitive and
therefore $\AA_2$ is transitive. To prove the strong closure of
$\AA_2$ is $\B(L^2(\M,\tau))$, by Theorem 3.1, we need  to prove
that $\AA_2$ is $2-$fold transitive. Let $Q\in Lat\AA_2^{(2)}$. Then
similar arguments as above show that $Q\in
Lat\AA_2^{(2)}=Lat\B(L^2(\M,\tau))^{(2)}$. By Proposition 1.2,
$\AA_2$ is $2-$fold transitive.
\end{proof}

\vskip 1.00cm

 {\bf Acknowledgements.}\,\, The authors thank the reviewers  for their valuable comments on this paper.


\vskip20pt

\begin{thebibliography}{99}
\bibitem{Ar} W.B. Arveson,   A density theorem for operator algebras,
\emph{Duke Math.J.} {\bf 34} (1967), 635-647.
\bibitem{Ar2} W.B. Arveson,  Subalgebras of $C\sp{*} $-algebras. II.
\emph{Acta Math.}  {bf 128}  (1972), no. 3-4, 271--308.
\bibitem{B-L} P. Biane, F.Lehner, Computation of some examples of
Brown's spectral measure in free probability, \emph{Colloq. Math.}
{\bf 90} (2001), 181-211.
\bibitem{Br} L.G. Brown,   Lidskii's theorem in the type II case,
Geometric methods in operator algebras, H. Araki and E. Effros
(Eds.) \emph{Pitman Res. notes in Math. Ser} {\bf 123}, Longman Sci.
Tech. (1986), 1-35.
\bibitem{D-H1} K. Dykema and U. Haagerup, Invariant subspaces of
Voiculescu's  circular operator, \emph{Geom. Funct. Anal.} {\bf 11}
(2001), 693-741.
\bibitem{D-H2} K. Dykema and U. Haagerup, Invariant subspaces of the
quasinilpotent DT-operator, \emph{J. Funct. Anal.} {\bf 209} no.2,
(2004), 332-366.
\bibitem{D-P} R. G. Douglas and C. Pearcy, Hyperinvariant subspaces and
transitive algebras.  \emph{Michigan Math. J.}  {\bf 19}  (1972),
1-12.
\bibitem{Fo} C. Foias, Invariant para-closed subspaces, \emph{Indiana
University Math. J.} {\bf 21} (1972), 887-906.
\bibitem{Fu} B. Fuglede,   A commutativity theorem for normal operators,
\emph{Proc. Nat. Acad. Sci. U.S.A.} {\bf 36} (1950), 35-40.
\bibitem{F-K} B. Fuglede and R.V.Kadsion, Determinant theory in finite
factors, \emph{Annals of Math.} {bf 55} (3)(1952), 520-530.
\bibitem{Ha1} U. Haagerup, The standard form of von Neumann algebras,
\emph{Math. Scand.} {\bf 37} (1975), 271-283.
\bibitem{Ha2} U. Haagerup, Spectral decomposition of all operators in a
$II_1$-factor, which is embeddable in $\R^{\omega}$. (Preliminary
version), MSRI 2001.
\bibitem{H-L} U. Haagerup and F. Larsen, Brown's spectral distribution
measure for $R-$diagonal elements in finite von Neumann algebras,
\emph{J.Funct.Anal.} {\bf 176} (2000), 331-367.
\bibitem{H-S} U. Haagerup and H. Schultz, Invariant Subspaces for
Operators in a General $II_1$-factor, preprint available at
http://www.arxiv.org/pdf/math.OA/0611256.
\bibitem{No} E. A. Nordgren, Transitive operator algebras.  \emph{J.
Math. Anal. Appl.}  {\bf 32}  (1970) 639--643.
\bibitem{NRR} E. A. Nordgren, H. Radjavi and P. Rosenthal, On density of
transitive algebras.  \emph{Acta Sci. Math. (Szeged)}  {\bf 30}
(1969), 175-179.
\bibitem{N-S} A. Nica and R. Speicher, $R-$diagonal pairs--a common
approach to Haar unitaries and circular elements. \emph{Fields
Institute Communications,} (1997), 149-188.
\bibitem{Ka} R.V. Kadison,  On the orthogonalization of operator
representation, \emph{Amer. J. Math.} {\bf 77} (1955), 600-621.
\bibitem{K-R} R. Kadison and J. Ringrose, ``Fundamentals of the Theory
of Operator algebras," Vols. 1, 2, Academic Press, INC, (1986).
\bibitem{Po} S. Popa, Orthogonal pairs of *-subalgebras in finite von
Neumann algebras, \emph{J. Oper. Theory,} {\bf 9} (1983), 253-268.
\bibitem{R-R} H. Radjavi and P. Rosenthal,  ``Invariant Subspaces'',
Springer-Verlag, New York, 1973.
\bibitem{S-S} P. Sniady and R. Speicher, Continuous family of invariant
subspaces for $R-$diagonal operators, \emph{Invent. Math.} {\bf 146}
(2001), 329-363.
\bibitem{VDN} D.V. Voiculescu, K. Dykema and A. Nica, ``Free Random
Variables", CRM Monograph Series, vol. 1, AMS, Providence, R.I.,
1992.
\end{thebibliography}
\end{document}